\documentclass{amsart}
\usepackage{amsfonts,amscd}
 
\newtheorem{claim}{}[section]
\newtheorem{theorem}[claim]{Theorem}
\newtheorem{lemma}[claim]{Lemma}
\newtheorem{proposition}[claim]{Proposition}
\newtheorem{corollary}[claim]{Corollary}

\def\proclaim #1. #2\par{\medbreak
\noindent{\bf#1.\enspace}{\sl#2}\par\medbreak}
\makeatother
 
\DeclareMathOperator{\ca}{C^*\text{-algebra}}

\DeclareMathOperator{\C}{{\mathcal C}}

\DeclareMathOperator{\cas}{C^*\text{-algebras}}

\begin{document}
 
\title{A double commutant theorem
for operator algebras}

\date{June 17, 2002.  Revised March 13, 2003.}

\author{David P. 
Blecher}\address{Department of Mathematics, University of Houston, Houston,
TX 77204-3008}
\email[David P. Blecher]{dblecher@math.uh.edu}

\author{Baruch Solel}\address{Department of Mathematics, Technion, 
32000 Haifa,
Israel}
\email[Baruch Solel]{mabaruch@techunix.technion.ac.il}  

\thanks{This research was supported in part by grants from
the National Science Foundation (Blecher) and the Fund for the
Promotion of Research at the Technion (Solel).}

\begin{abstract}
Every unital nonselfadjoint operator 
algebra possesses 
canonical and functorial classes of faithful 
(even completely isometric) Hilbert 
space representations satisfying 
a double commutant theorem generalizing von Neumann's 
classical result. 
Examples and complementary results are given.
\end{abstract}
 
\maketitle

\let\text=\mbox

\section{Introduction.}
In this paper we consider possibly nonselfadjoint norm closed algebras
of operators
on a Hilbert space $H$.
The {\em general theory} of such {\em operator algebras}, and
of their representations, is rather sparse in contrast to the selfadjoint 
case, namely
the $C^*$-algebra theory.  The contrast is easily seen
in the lack of certain fundamental tools which are
available in the selfadjoint case, such as
von Neumann's double commutant theorem or Kaplansky's 
density theorem.  In recent years, with the
help of operator space theory, the situation has changed somewhat,
and the general theory  of operator algebras has been growing
rapidly.  The theory that is emerging also has the feature that it
links much more closely together the three subjects of operator algebras,
$C^*$-algebra theory, and the theory of rings and modules.
The present paper is an attempt in this spirit to establish a
double commutant theorem for general operator algebras, one that
would be useful at least for tackling certain problems.
Our main result simultaneously resembles two classical and 
fundamental results: von Neumann's double commutant theorem
\cite{vN},
and Nesbitt and Thrall's purely algebraic result \cite{NT}
 to the effect that any module $M$ over a ring $R$, which is
a `generator for $_RMOD$',
satisfies the appropriate double commutant theorem.
Although for a 
subalgebra $A \subset B(H)$ 
the double commutant $A''$ may not
equal the weak* closure $\overline{A}^{weak*}$
within $B(H)$,
we show that 
important and canonical classes
 of completely isometric
Hilbert space representations $\pi : A \rightarrow B(K)$ of $A$  do have
the {\em double commutant
property}, by which we will mean that
$\pi(A)'' = \overline{\pi(A)}^{weak*}$ in $B(K)$.   
In fact for a certain subclass of these representations
this double commutant is also  
isomorphic to $A^{**}$.
If $A$ is a `dual operator algebra', then we can find 
canonical classes of completely isometric normal representations
with the property $\pi(A)'' = \pi(A)$, as one is accustomed
to for von Neumann algebras.   

A good illustration is the very  simple example
${\mathcal T}_2$, the algebra of $2 \times 2$ upper triangular
matrices.  The usual representation of ${\mathcal T}_2$
 on $\mathbb{C}^2$
clearly  does not have the double commutant property.
However the direct sum $\rho$ of the usual representation and the
1-dimensional representation consisting of evaluation at
the 1-1 entry,
does possess the double commutant property.
Indeed this representation $\rho$ satisfies our general `sufficient
condition' for the double commutant property (Theorem \ref{pThm}), but
the usual representation of ${\mathcal T}_2$ does not.

The substitution
of a given embedding $A \subset B(H)$ of an operator algebra $A$
for another in which the
double commutant theorem holds, is easily justified by 
the trend in the recent
general theory of operator algebras towards a `coordinate-free' approach,
i.e. not to focus on any
one fixed embedding $A \subset B(H)$.  In the new perspective
one is encouraged to think about all representations of $A$ 
simultaneously; some may be better than others for solving 
certain problems.

Before we turn to specific details, we feel it is important to mention that
although some of the representations we consider are 
quite `large', in 
some applications
this should not matter - the important thing is often just
the universal property
and the functoriality, and the strong links to the $C^*$-algebra theory.
A good illustration of this may be found in
the paper \cite{BFTM}, where a difficult
problem concerning nonselfadjoint operator algebras was solved by
transferring it to the $C^*$-algebra world using some rather large 
representations.  We will give other such
applications of our results in a sequel paper, as well as a
characterization of the double commutant property which is
quite different to the considerations in the present paper.  

We now turn to specific details, and recall some definitions
and facts (see e.g. \cite{BLM} for more details if needed).
Although everything in this
paper can be done within the Banach algebra context
and without the matrix norms of operator space theory,
for specificity we use the operator space 
context.   Thus following the lead suggested by Arveson 
\cite{Arv1}, and by operator space theory,
we define an {\em abstract operator algebra} 
- or simply an {\em operator algebra} -  to be a
Banach algebra $A$ which is also an operator 
space, such that there exists a completely 
isometric homomorphism $\theta : A \rightarrow B(H)$ 
for some Hilbert space $H$.   There is an abstract characterization
of these operator algebras due to the first author with Ruan and 
Sinclair, but we shall not need this here.   Except in the last section,
operator algebras are assumed to
have a contractive approximate identity (c.a.i.).  If $A$ has an
identity of norm 1 then we say that $A$ is {\em unital}.

It is helpful to use the language of Hilbert 
modules (see for example 
\cite{Ri,DP,MSmem}).  
For the purposes of this paper, we define a Hilbert
$A$-module to be  a Hilbert space $H$ which is a
nondegenerate $A$-module via a 
completely contractive nondegenerate 
representation $\pi : A  \rightarrow B(H)$.
Thus $a \zeta = \pi(a)(\zeta)$, for $a \in A, \zeta \in H$.  
The theory below may be modified to include 
`contractive Hilbert $A$-modules', we omit the easy details.
Except in the final 
section of our paper,
the word `nondegenerate' used above means that 
the span of such products $a \zeta $ is dense in $H$.
This correspondence between Hilbert $A$-modules and 
completely contractive nondegenerate representations is bijective.
Henceforth we use the term `representation of $A$' for
 a completely contractive nondegenerate representation.
We write $_AHMOD$ for the category of Hilbert $A$-modules,
with morphisms the bounded maps which intertwine the 
representations (that is, the bounded $A$-module maps).  
This category is closed under direct sums and quotients
by closed submodules (see e.g. \cite{MSmem}).
If $\alpha$ is
a cardinal then the direct sum of $\alpha$ copies of a
Hilbert $A$-module $H$, or of its associated representation
$\pi : A  \rightarrow B(H)$, is called a {\em multiple} of
$H$ or of $\pi$; and is written as $H^{(\alpha)}$ or $\pi^\alpha$.     
We say that two Hilbert $A$-modules are {\em spatially equivalent},
and write $H \cong K$, if they are isometrically $A$-isomorphic,
that is if there exists a unitary $A$-module map from 
$H$ onto $K$.  We say that a closed $A$-submodule $K$ of an
Hilbert $A$-module $H$ is {\em $A$-complemented} if the 
projection of $H$ onto $K$ is an $A$-module map.  This may be 
reformulated in several equivalent ways (see for example 
the discussion in \cite{MSmem}).   We  say that representations 
$\pi , \theta$ of $A$ are
{\em quasi-equivalent} if there is a multiple of $\pi$ which is
spatially equivalent to a multiple of $\theta$.   Thus
two Hilbert $A$-modules are  quasi-equivalent if and only if
there are cardinals $\alpha$ and $\beta$ (which we may clearly
assume to be equal) such that $H^{(\alpha)} \cong K^{(\beta)}$.   
 We say that a
Hilbert $A$-module $H$ is {\em $A$-universal},
if every $K \in \; _AHMOD$ is
isometrically $A-$isomorphic (that is,
spatially equivalent) to an $A$-complemented
submodule of a direct sum of copies of $H$.
We  say that a module $H \in \;
_AHMOD$ is a {\em generator} (resp.  {\em cogenerator})
for $_AHMOD$ if for every nonzero morphism $R : K \rightarrow L$
of $_AHMOD$, there
exists a  morphism $T : H  \rightarrow K$
(resp. $T : L  \rightarrow H$) of $_AHMOD$ with $R T \neq 0$
(resp. $T R  \neq 0$).
We will say that $H$ is {\em sub-tracing}
if the definition above for
generator  is modified so that
$K$ ranges over the set of submodules of $H$.  
We say that  $H$ is {\em completely sub-tracing} if 
a countably infinite multiple of $H$ is sub-tracing.  
We also use these
terms when referring to the associated representation on $H$.
Thus for example we will often refer to a {\em representation}
 $\pi :
A  \rightarrow B(H)$ as being $A$-universal, or sub-tracing.
One would expect, just as in pure algebra,
that there are many useful alternative characterizations
of generators, cogenerators, and
sub-tracing modules.  Indeed we will provide some 
in Section 2.   

Note that any generator is sub-tracing.  Since any multiple of a
generator is also a generator, we see that 
any generator is completely sub-tracing.
Also, any $A$-universal Hilbert module $H$ is a generator.  To see this,
suppose that $T : K \rightarrow L$
is a nonzero $A$-module map.   W.l.o.g., there is a
cardinal $\alpha$ such that $K$  is an $A$-complemented
submodule of $H^{(\alpha)}$;  let $Q$ be the
associated projection onto $K$ from  $H^{(\alpha)}$.
Let $\epsilon_i$ be
the inclusion map of $H$ into $H^{(\alpha)}$  as
its $i$th summand $H$.  If every map
$T \circ Q \circ \epsilon_i$ is
zero, then $T = 0$, which is a contradiction.
 
The main results of our paper are the following:

\begin{theorem}  \label{pThm}  Let $A$ be an operator
algebra with a contractive approximate identity.
A Hilbert $A$-module which is
a generator or cogenerator for 
$_AHMOD$, or which is completely sub-tracing, has the 
double commutant property.
\end{theorem} 

\begin{theorem}  \label{Thm}  Let $A$ be an operator 
algebra with contractive approximate identity. 
\begin{itemize}
\item [(1)]  There exist $A$-universal representations
for $A$.
\item [(2)]  Any two $A$-universal representations 
for $A$ are quasi-equivalent.
 \item [(3)]  If $\pi$ is a 
representation of $A$ which is 
quasi-equivalent to an $A$-universal representation,
then $\pi$ is an $A$-universal representation.
\item [(4)]  If $\pi$ is an $A$-universal representation
of $A$ on a Hilbert space $H$, then 
$$\pi(A)'' = \overline{\pi(A)}^{weak*}.$$
\item [(5)]   If $\pi$ is an $A$-universal representation,
 then $A^{**}$ is isomorphic to $\overline{\pi(A)}^{weak*}$
 via a completely isometric weak*-homeomorphic homomorphism
$\rho : A^{**}  \rightarrow \overline{\pi(A)}^{weak*}$
such that $\rho(\hat{a}) = \pi(a)$ for all $a \in A$.
\end{itemize}
\end{theorem}

\begin{proof}  We prove only items (2)-(4) now, 
and defer the proofs of the other assertions.  In fact
(3) is clear by the definitions, or is an easy exercise.
The proof of (2) is a simple application of set theory,
and the well known `Eilenberg swindle'.
If $H$ and $K$ are two $A$-universal representations,
then there exist cardinals $\alpha$ and $\beta$, and 
Hilbert $A$-modules $M$ and $N$, such that 
$H \oplus M \cong K^{(\alpha)}$ and $K \oplus N 
\cong H^{(\beta)}$.    Without loss of generality,
by adding on extra multiples of $H$ or $K$ to the last 
two equations, $\alpha = \beta$.  We may also assume that 
$\alpha$ is a large enough cardinal so that 
$\alpha \cdot \alpha$ equals $\alpha$.
Then 
$$K^{(\alpha)} \cong K^{(\alpha)} \oplus K^{(\alpha)}
\oplus \cdots
\cong H \oplus M \oplus H \oplus M \oplus \cdots .$$ 
By associativity we get 
$$K^{(\alpha)} \cong H \oplus K^{(\alpha)} \oplus K^{(\alpha)} 
\oplus \cdots \cong H \oplus K^{(\alpha)}.$$
Since $\alpha \cdot \alpha = \alpha$, a multiple of the 
last equation yields 
$$K^{(\alpha)} \cong H^{(\alpha)} \oplus K^{(\alpha)}.$$
Similarly, $H^{(\alpha)} \cong H^{(\alpha)} \oplus K^{(\alpha)}
\cong K^{(\alpha)},$ which proves (2).

Item (4) follows from Theorem \ref{pThm} and the fact 
above
that any $A$-universal Hilbert module $H$ is a generator.  
\end{proof}

It follows from (4) and (5) that for $A$-universal representations,
there is an automatic `Kaplansky density' result,
which is really Goldstine's lemma in disguise.
 
The main results above are proved in the first few sections
of our paper.
In Sections 4--6 we give examples and complementary results.
For example we study there dual algebras;
the relations between various classes of
Hilbert modules; and in the 
final section we study 
operator algebras without c.a.i., 
showing for example that all
of Theorem \ref{Thm} with the exception of (4) holds in
complete generality.  

We list now some background facts that we will make 
much use of (often without 
comment).  One fact which is of great
assistance when dealing with operator 
algebras with c.a.i. but no identity,  is the following:
if $B$ is a $C^*$-algebra generated (as 
a $C^*$-algebra) by a closed subalgebra $A$ which has a c.a.i., 
then 
any $b \in B$ is a
product $a b'$ (or $b' a$) with $a \in A, b' \in B$.   Equivalently:
\begin{equation} \label{anc}
\text{Any c.a.i. for} \; A \; \text{
is also one for the} \;
\ca \; \text{generated by A} 
\end{equation}
See e.g. \cite{BLM} Chapter 2.   We use 
the notation $[A K]$ for the norm closure of the 
span of products of a term in $A$ with a term in $K$.     

If ${\mathcal S} \subset B(H)$ then we define
${\mathcal S}^* = \{ x^* : x \in {\mathcal S} \}$.  If 
$K$ is another Hilbert space  (resp. if $\gamma$ is a cardinal),
then we write ${\mathcal S} \otimes I =
 \{ x \otimes I :  x \in {\mathcal S} \}$ for the 
set of appropriate 
`multiples' of elements in ${\mathcal S}$.   This is
a set of operators on $H \otimes K$ (resp.  on $H^{(\gamma)}$).
It is a simple computation that the following relations hold: 
$$\overline{{\mathcal S} \otimes I}^{weak*} 
= \overline{{\mathcal S}}^{weak*} \otimes I,$$
$$\left( {\mathcal S} \otimes I \right)'' =  
{\mathcal S}''  \otimes I,$$
 $$\overline{{\mathcal S}^*}^{weak*}
= \left( \overline{{\mathcal S}}^{weak*} \right)^*,$$
and 
$$\left( {\mathcal S}^* \right)'' =
\left( {\mathcal S}'' \right)^*. $$
Hence ${\mathcal S}$ has the double commutant property if and
only if ${\mathcal S}^* $ has the double commutant property,
and if and only if ${\mathcal S} \otimes I$ 
has the double commutant property.

\section{Generators and traces}

Following algebra texts (e.g. \cite{AF} p. 109)
if $H, K$ are Hilbert $A$-modules, then we define the {\em trace}
$Tr_K(H)$ to be the
closure of the set of finite sums of elements taken
from the ranges of bounded $A$-module
maps from $H$ into $K$.   We define the {\em 
reject}  $Rej_K(H)$ to be the intersection of the kernels of 
all bounded $A$-module maps from  $K$ into $H$.
Clearly $Tr_K(H)$ and $Rej_K(H)$ are 
closed submodules of $K$.  

\begin{lemma}  \label{tr}  Let $H$ be a 
Hilbert $A$-module.  Then 
\begin{itemize}
\item [(1)]   $H$ is a generator
for $_AHMOD$  if and only if $Tr_K(H) = K$ for all
Hilbert $A$-modules $K$.
\item [(2)]  $H$ is a cogenerator
for $_AHMOD$  if and only if $Rej_K(H) = \{ 0 \}$
 for all
Hilbert $A$-modules $K$.
\item [(3)]  $H$ is sub-tracing  if and only if $Tr_K(H) = K$ for  
closed submodules $K$ of $H$.
\end{itemize}  
  \end{lemma}
 
\begin{proof}   (1).   Suppose that $H$ is a generator.
If $Tr_K(H) \neq K$, then there exists a
nonzero bounded $A$-module map $R : K \rightarrow K/Tr_K(H)$
annihilating $Tr_K(H)$, namely the  quotient map. 
   Since the quotient of
Hilbert $A$-modules is a Hilbert $A$-module, and since
$H$ is a generator, there exists a $T \in B_A(H,K)$
with $R T \neq 0$.  This contradicts the definition of the trace.
(3) is proved similarly to (1).
 
Conversely,  suppose that $Tr_K(H) = K$, and
$R : K \rightarrow L$ is a bounded $A$-module map.
If $R \circ T = 0$ for all $T \in B_A(H,K)$ then  $R$ is
zero on  $Tr_K(H) = K$.  So $H$ is a generator.

For (2), assume that $Rej_K(H) = \{ 0 \}$.
If $R : L \rightarrow K$ is a nonzero morphism,
but that $T R = 0$ for all $T \in \; _AB(K,H)$, then
$R$ maps into $Rej_K(H) = \{ 0 \}$.  Hence $R = 0$.
Conversely, if $Rej_K(H) \neq \{ 0 \}$, then the 
inclusion map $\epsilon : Rej_K(H) \rightarrow K$ is a
nonzero bounded $A$-module map with $T \epsilon = 0$
for all $T \in \; _AB(K,H)$.
\end{proof}

\begin{lemma} \label{bfge}  Suppose that $\pi$
is a generator for $_AHMOD$, that $\sigma$ is a
cogenerator for $_AHMOD$, and  that
$\rho$ is any representation in $_AHMOD$.  We have:
\begin{itemize}
\item [(1)]  $\pi \oplus \rho$ is  a generator for $_AHMOD$,
$\sigma  \oplus \rho$ is  a cogenerator for $_AHMOD$,
and $\pi \oplus \sigma 
\oplus \rho$ is both a
generator and cogenerator for $_AHMOD$.
\item [(2)]  $\pi$ and $\sigma$ are 1-1 (i.e. faithful).
\item [(3)]   $\rho$ is a
generator (resp. cogenerator) for $_AHMOD$ if and only if 
a multiple of $\rho$ is a
generator (resp. cogenerator) for $_AHMOD$.
\item [(4)]  If $\rho$ is quasi-equivalent to 
$\pi$ (resp. to $\sigma$) then 
$\rho$ is a generator  (resp. cogenerator)
for $_AHMOD$.
\end{itemize}
  \end{lemma}
 
\begin{proof}   (1) is obvious from the definitions.
For (2), suppose that  $\pi(a) = 0$, and let 
$K$ be a faithful $A$-module, with corresponding 
1-1 representation $\sigma$.  If $T \in \; _AB(H,K)$
then $\sigma(a) T(\zeta) = T \pi(a) \zeta = 0$.  Hence
$\sigma(a)$ is zero on $Tr_K(H) = K$ (using 
Lemma \ref{tr}).  So $a = 0$.   A similar obvious 
argument proves the assertion for $\sigma$.   
We leave (3) and (4) as exercises; they will not be 
explicitly used in the paper.
\end{proof}    

\begin{proposition} \label{Prop1}   Let $A$
be an operator algebra with c.a.i..  Suppose that $\rho : A \rightarrow B(H)$
is a sub-tracing representation.  Then $\rho(A)'' \subset \;
\text{alg lat } \; \rho(A)$.
\end{proposition}
 
\begin{proof}   
Fix $x \in H, x \neq 0$,
and consider the closed span $K$ of
$\rho(A) x$ in $H$.   
If $\{ e_\alpha \}$ is a c.a.i. for $A$ then
$\rho(e_\alpha) x \rightarrow x$.  Thus $x \in K$.
Suppose that $T \in \rho(A)''$.
If $V \in \; _AB(H,K)$ then $V$ (regarded as a map into $H$)
is in $\rho(A)'$.  Thus $T V H = V T H \subset K$.  Consequently,
$T$ maps the trace $Tr_K(H)$ into $K$.   By Lemma \ref{tr} (3),
$T(x) \in T(K) \subset K$.  Since 
$x$ was arbitary, we are done.
\end{proof}

\begin{corollary} \label{dct1}  If $A$ is an operator algebra with c.a.i.
and if $\pi$ is a completely sub-tracing representation of
$A$ on a Hilbert space $H$, then
$$\pi(A)'' = \overline{\pi(A)}^{weak*}.$$
\end{corollary}
 
\begin{proof}    By definition, for a separable infinite 
dimensional Hilbert space $H_0$ we have that 
$\pi(\cdot) \otimes I_{H_0}$ is sub-tracing.
Hence by the previous proposition we have
$$(\pi(A) \otimes I_{H_0})'' \subset \; \text{alg lat }
\left( \pi(A) \otimes I_{H_0} \right)
 \subset \; \text{alg lat } \left( \overline{\pi(A)}^{weak*}
 \otimes I_{H_0} \right).$$
By well known facts about `reflexive algebras' 
(see e.g. Lemma 15.4 in \cite{DNA}),
 $$\text{alg lat } \left( \overline{\pi(A)}^{weak*}
 \otimes I_{H_0} \right) =
\overline{\pi(A)}^{weak*} \; \otimes \; I_{H_0}.$$
Of course $(\pi(A) \otimes I_{H_0})'' =
\pi(A)''  \otimes I_{H_0}$.  Putting the facts above together
yields
$$\pi(A)''  \otimes I_{H_0}  \subset \overline{\pi(A)}^{weak*} \;
\otimes \; I_{H_0}$$
so that
$$\pi(A)''   \subset \overline{\pi(A)}^{weak*}.$$
The other direction is trivial since
$\pi(A)''$ is weak* closed.
\end{proof}
 
By the last result,
we are now almost done with the 
proof of  Theorem \ref{pThm}.   The final part is completed 
as follows.  Suppose that $H$ is a cogenerator for 
$_AHMOD$.  Then by 
simple observations in the next section
(before (\ref{adj})), $H$ is a generator for 
$_{A^*}HMOD$.  Thus the image of $A^*$ in $B(H)$ has the 
double commutant property.  The facts at the 
end of Section 1 now complete the proof.

\vspace{3 mm}

{\bf Remark.}  It is fairly clear that the definitions of 
$A$-universal, generator, cogenerator, or sub-tracing, are
{\em functorial}.  In 
particular, if $B$ is another such operator algebra,
and if the categories $_AHMOD$ and $_BHMOD$ are equivalent as
categories, then it is easy algebra to check that 
the equivalence functor takes $A$-universal 
representations to $B$-universal representations and vice versa.
Similarly for generators, cogenerators, or sub-tracing
Hilbert modules.

\section{The universal $C^*$-algebra and
universal representations}

We will need to recall several simple facts
(see e.g. \cite{BLM} for more details, and examples, if needed).
Firstly, there is a canonical functor $A \mapsto
C^*(A)$ from the category of 
operator algebras (with c.a.i.) and completely contractive
homomorphisms, to the category of $C^*$-algebras
and *-homomorphisms, with the following universal property:
there exists a completely
isometric homomorphism $i : A \rightarrow C^*(A)$ such that $i(A)$
generates $C^*(A)$ as a $C^*$-algebra, and such that if $\phi : A
\rightarrow D$ is any completely contractive homomorphism into a
$C^*$-algebra $D$, then there exists a (necessarily unique)
*-homomorphism $\tilde{\phi} : C^*(A) \rightarrow D$ such that
$\tilde{\phi} \circ i = \phi$.  The algebra $C^*(A)$ is called
the {\em maximal $C^*$-algebra} generated by $A$, and is
sometimes written as $C^*_{max}(A)$.  
For those interested in algebra,
this universal property essentially says that the
functor  $A \mapsto C^*(A)$ is the left adjoint to the forgetful
functor from the category of $C^*$-algebras to 
the category of operator algebras.

From the universal property of
$C^*(A)$, it is clear that if $H$ is a
Hilbert $A$-module, then the associated representation
$\pi$ has a unique extension $\tilde{\pi}$ which is
a nondegenerate *-representation of $C^*(A)$ on $H$.
Conversely every nondegenerate *-representation of $C^*(A)$ on $H$
restricts (using the fact (\ref{anc}) from
Section 1 if necessary) to a nondegenerate representation of $A$ on $H$.  
Thus we may regard Hilbert $A$-modules as Hilbert $C^*(A)$-modules,
and vice versa, in this canonical way.
By symmetry every Hilbert $A$-module is
also a nondegenerate Hilbert module over the subalgebra $A^*$ of
$C^*(A)$ (one may deduce the nondegeneracy from 
fact (\ref{anc}) from Section 1 again).   
If $T : H \rightarrow K$ is  a bounded  $A$-module map,
then it is easy to see that
$T^* : K \rightarrow H$ is an $A^*$-module map; and conversely.
From this we can make a few simple deductions.  Firstly,
it follows from the last fact that $H$ is a generator for 
$_AHMOD$ if and only if $H$  is a cogenerator for
$_{A^*}HMOD$.    Similarly,  $H$  is a cogenerator for
$_AHMOD$ if and only if $H$  is a generator for
$_{A^*}HMOD$.    
Secondly, if we  call a bounded $A$-module
map  $T$ between Hilbert $A$-modules {\em adjointable} if
$T^*$ is also an $A$-module map, then we have from the above that: 
\begin{equation} \label{adj}
T \; \text{is adjointable if and only if} \; T \; \text{is a}
\; C^*(A)\text{-module map.} 
\end{equation} 
We shall not need adjointable maps very much,
except in the special case
that $i$ is an isometric $A$-module map between
Hilbert $A$-modules, such that $i^*$ is an $A$-module map.
It follows from (\ref{adj}) that
$i$ and $i^*$ are $C^*(A)$-module maps.  In particular 
we deduce that {\em unitary}
morphisms, i.e. unitary $A$-module maps, are $C^*(A)$-module maps.
That is, two Hilbert $A$-modules are spatially equivalent as
$A$-modules if and only if they are spatially equivalent as
$C^*(A)$-modules.

From the above, it also
follows that the class of $A$-complemented submodules of a
Hilbert $A$-module $H$ is the same as the class of
closed $C^*(A)$-submodules of $H$.   Hence  Hilbert $A$-module direct sums
(resp. summands) of Hilbert $A$-modules
are the same as Hilbert $C^*(A)$-module
direct sums (resp. summands).   From these considerations the following
result is
clear.  Note that part (2) below  establishes (1) of Theorem \ref{Thm}: 

\begin{corollary} \label{uni}   
Let $A$
be an operator algebra with c.a.i.. 
\begin{itemize}
\item [(1)]   A Hilbert 
$A$-module $H$ is $A$-universal if and only if
$H$ is $C^*(A)$-universal.  
\item [(2)]  Any
$C^*(A)$-universal representation  of $C^*(A)$, such as the 
usual `universal representation' $\pi_u$ of $C^*(A)$,
restricts to an $A$-universal representation of $A$.
\item [(3)]  Two Hilbert $A$-modules are quasi-equivalent as
Hilbert $A$-modules if and only if they are quasi-equivalent as
Hilbert $C^*(A)$-modules.
\item [(4)]   If $\pi$ and $\theta$ are 
quasi-equivalent representations of $A$, then
there exists a (necessarily unique)
weak*-homeomorphic completely isometric isomorphism
$\rho : \overline{\pi(A)}^{weak*} \rightarrow 
\overline{\theta(A)}^{weak*}$
such that $\rho \circ \pi = \theta$. \end{itemize}   
\end{corollary}

\begin{proof} (1)-(3) are obvious from the discussion above.
Let $\C = C^*(A)$.
If $\pi$ and $\theta$ satisfy the hypothesis in (4),
then by (3) we know that $\tilde{\pi}$ and $\tilde{\theta}$
are quasi-equivalent in the usual $C^*$-algebraic sense.  
Thus by  5.3.1 (ii) in \cite{Dix}, there
is a $W^*$-isomorphism $\Phi : \overline{\tilde{\pi}(\C)}^{weak*} \rightarrow
\overline{\tilde{\theta}(\C)}^{weak*}$ such that
$\Phi \circ \pi = \theta$.  The restriction of $\Phi$ to
$\pi(A)$ maps onto $\theta(A)$, so that 
by weak*-continuity we obtain the result.  
\end{proof}

Thus there is a `canonical' $A$-universal representation 
of $A$, namely the restriction of the 
universal representation  of $C^*(A)$ to $A$.  We will 
call this {\em the  universal representation of $A$}, and 
we write this representation of $A$ as $\pi_u$.  
  
The following result, which shall not be used in an essential
way in this
paper, shows that the 
universal representation satisfies
 quite a strong form of the 
`sub-tracing' condition.
 
\begin{proposition} \label{Prop1b}   Let $A$
be an operator algebra with c.a.i..  Suppose that $\rho :
C^*(A) \rightarrow B(H)$ is a nondegenerate
*-representation with the following
property: For every state $\varphi$ on $C^*(A)$ there exists
a $\xi \in H$ such that $\varphi(b) = \langle \rho(b) \xi , \xi \rangle$
for all $b \in C^*(A)$.  Then for every topologically singly 
generated $A$-submodule 
$K$ of $H$, there is a partial isometry in 
$A'$ with range $K$.
\end{proposition}

\begin{proof}   Let $B = C^*(A)$.
Fix $x \in H, \Vert x \Vert = 1$,
and consider the closed span $K$ of
$\rho(A) x$ in $H$.   As noted in the proof of 
\ref{Prop1}, $x \in K$.
Since $\rho(A) K  \subset K$, by the universal
property of $C^*(A)$ there exists a *-representation $\pi$ of
$C^*(A)$ on $K$ with $\pi(a) = \rho(a)_{|_K}$ for $a \in A$.
Let $\varphi = \langle \pi(\cdot) x , x \rangle$.   This
is a state, so by hypothesis there exists a  $\xi \in H$ such that
$\varphi = \langle \rho(\cdot) \xi , \xi \rangle$
for all $b \in B$.   Therefore
$$\Vert \rho(b) \xi \Vert^2 = \Vert \pi(b) x \Vert^2$$
for all $b \in B$, and so there is a well defined unitary
$V_0$ from $[\rho(B) \xi]$ to $[\pi(B) x] \subset K$ taking
$\rho(b) \xi$ to $\pi(b) x$.  Extend $V_0$ to a partial isometry
$V \in B(H)$ by setting it to be zero on $[\rho(B) \xi]^\perp$.
It is easy to see that
$$\rho(A) [\rho(B) \xi]^\perp
\subset \rho(B) [\rho(B) \xi]^\perp \subset [\rho(B) \xi]^\perp,$$
and therefore $V \rho(a) y = \rho(a) V y = 0$ if $y \in [\rho(B) \xi]^\perp$
and $a \in A$.  
On the other hand, if $y = \rho(b) \xi \in [\rho(B) \xi]$, 
then $$V \rho(a) y = V \rho(ab) \xi = \pi(ab) x
= \rho(a) \pi(b) x =  \rho(a) V \rho(b) \xi = \rho(a) V y.$$
Hence $V \in \rho(A)'$. 
\end{proof}     

{\bf Remark:}  If $\pi_u$ is the universal representation of
$C^*(A)$ on $H_u$, and if $H_0$ is any Hilbert space, then 
the *-representation $\rho = \pi_u(\cdot) \otimes I_{H_0}$
of $C^*(A)$ on $H_u \otimes H_0$, satisfies the conditions
of the proposition.

\section{Dual operator algebras and normal representations}

In this section we 
turn to dual operator algebras, and we 
will
also prove the remaining part, namely (5),
of Theorem \ref{Thm}.
We again begin by recalling some facts
and notations (see e.g. \cite{BLM} for more details 
if needed).
A {\em dual operator algebra} is an
operator algebra $A$ which has a predual such that
$A$ is completely isometrically isomorphic, via a
 homomorphism which is a homeomorphism for the weak*
topologies, to a
$\sigma$-weakly closed unital subalgebra of $B(H)$.
There is an abstract characterization of dual operator algebras
due to Le Merdy, with a contribution by the first author,
but we shall not need this here.
A {\em normal representation} of a dual operator algebra
is a unital completely contractive weak* continuous
homomorphism $\pi :  A \rightarrow B(K)$.  
We write $_A{NHMOD}$ for the category of the
Hilbert modules corresponding to such
normal representations, and call an object in
$_A{NHMOD}$ a {\em normal Hilbert $A$-module}.  The morphisms
are the same as in $_AHMOD$.   Again it is a simple exercise that
$_A{NHMOD}$ is closed under direct sums.

Let  $A$ be an operator algebra with c.a.i..
It is a well known fact (that appears first in
\cite{ER}, and which may be deduced 
for example from the first part of the next proof) 
that $A^{**}$ is a unital 
operator algebra in a canonical way. 
For any $H \in \; _AHMOD$, with corresponding
representation $\pi$, we may use the universal
property of $C^*(A)$ to get a *-representation
$\tilde{\pi} : C^*(A) \rightarrow B(H)$.   As is explained
in any text on $C^*$-algebras, we may extend
this *-representation in a unique fashion
to a unital normal *-representation
$\theta : C^*(A)^{**} \rightarrow B(H)$.
Let $\bar{\pi}$ be
$\theta$ restricted to $A^{**}$.  Then clearly
$\bar{\pi}$ is the unique normal representation
$A^{**} \rightarrow B(H)$ extending $\pi$ from
$A$.   Moreover, since
$\bar{\pi}$ is w*-continuous, its range is contained in
the dual operator algebra $\overline{\pi(A)}^{weak*}$.

We remark in passing that the converse is true too:
any
normal representation $\rho : A^{**} \rightarrow B(H)$
restricts to a
nondegenerate completely
contractive representation of $A$ on $H$.
In fact, generalizing
a well known fact for $C^*$-algebras (see p. 53
of \cite{Ri}), it is clear that the categories $_AHMOD$ and
$_{A^{**}}NHMOD$ are equivalent for any
operator algebra  $A$ with c.a.i..
 
\begin{proof}  {\em (of (5) of Theorem \ref{Thm}.)}
Applying the remark above to the
universal representation $\pi_u$ of $A$, we obtain a
normal representation $\overline{\pi_u} : A^{**}
\rightarrow \overline{\pi_u(A)}^{weak*} \subset B(H_u)$.
In fact $\overline{\pi_u}$ is completely isometric, since it
is the restriction of the faithful *-isomorphism between
$\C^{**}$ and $\pi_u(\C)''$, where $\C =
C^*(A)$.  Thus by the
Krein-Smulian theorem, the image $B$ of $A^{**}$ under
$\overline{\pi_u} $ is weak* closed, and
$\overline{\pi_u}$ is a homeomorphism for the weak* topologies.
Since $B$ contains $\pi_u(A)$ we have that
$\overline{\pi_u(A)}^{weak*} = B$.
This proves the result in the
special case that $\pi = \pi_u$.

To prove the general case, suppose
that $\theta$ is an $A$-universal representation of $A$.
By (2) of Theorem \ref{Thm}, $\theta$ is quasi-equivalent to
$\pi_u$.  By (4) of \ref{uni}, there exists a
weak* homeomorphic completely isometric isomorphism
$\rho : \overline{\pi_u(A)}^{weak*} \rightarrow
\overline{\theta(A)}^{weak*}$
such that $\rho \circ \pi = \theta$.    Composing $\rho$ with
the map $\overline{\pi_u} : A^{**} \rightarrow \overline{\pi_u(A)}^{weak*}$
of  the previous paragraph, gives the desired map in (5)
of Theorem \ref{Thm}.
\end{proof}

As a corollary of (4) and (5) of Theorem \ref{Thm},
one may immediately obtain the following fact  
which implies some results proved in \cite{BOMD} and \cite{LMsa}:

\begin{corollary}  Suppose that 
$A$ is an operator algebra with c.a.i., which possesses
an $A$-universal representation $\pi$ with $\pi(A)'$ 
selfadjoint.   Then $A$ is a $C^*$-algebra.
\end{corollary}

\begin{proof}   If $\pi(A)'$ is selfadjoint,
then by the above $A^{**}  \cong \pi(A)''$ is a $W^*$-algebra.
The proof is completed by an appeal to the following Lemma.
\end{proof}  

\begin{lemma} \label{acs}  Suppose that $A$ is an operator algebra such that
$A^{**}$ possesses an involution with 
respect to which $A^{**}$ is a $C^*$-algebra.  Then $A$ is a $C^*$-algebra.
\end{lemma}

\begin{proof}   Suppose that $A$ is a subalgebra of
a $C^*$-algebra $B$.
Then we have the following closed subalgebras:
$\hat{A} \subset A^{**} \subset B^{**}$.   It is well known 
that a contractive homomorphism between $\cas $ is a 
*-homomorphism. 
Thus, if $A^{**}$ possesses an involution as stated, then 
it follows that $A^{**}$ is a *-subalgebra of $B^{**}$.  We 
complete the proof by showing that 
$\hat{A}$ is closed under the above involution.
For if $a \in A$, then  $\hat{a}^*
\in A^{**} \cap \hat{B}$.   By a basic fact 
for Banach spaces, $A^{**} \cap \hat{B} = \hat{A}$.
\end{proof}

A different proof of this lemma was found together with Le Merdy
around '99.

If $\pi$ is a normal completely isometric representation
of a dual operator algebra,
then it follows from the Krein-Smulian theorem that 
$\pi(A)$ is weak* closed.   If further $\pi$  
is completely sub-tracing, then it follows immediately 
from  Corollary \ref{dct1} that 
$$\pi(A)'' = \pi(A).$$
One may define {\em normal generators} and 
{\em normal $A$-universal representations} for the 
category
$_ANHMOD$, in an obvious way.  It is clear as before that every
normal $A$-universal is a 
normal generator, and that every 
normal generator is completely sub-tracing.    

In order to see that for any dual operator algebra
there do exist normal $A$-universal representations, 
and for its own intrinsic interest, 
we will define a
{\em maximal $W^*$-algebra} $W^*(A)$ of a  dual operator algebra
$A$.    This is a
$W^*$-algebra, together with a weak* continuous completely isometric
homomorphism $j : A \rightarrow W^*(A)$ whose range
generates $W^*(A)$ as a $W^*$-algebra,  
and which possesses the
following universal property:
given any normal representation
$\pi : A \rightarrow B(H)$, there exists a
normal *-representation $W^*(A) \rightarrow B(H)$
extending $\pi$. 
It is elementary to define this
if $A = B^{**}$ for an operator algebra $B$, in this case simply
let $W^*(A) = C^*(B)^{**}$, and one may easily check that
this has the desired universal property.  But if $A$ is a
general dual operator algebra a little more care is 
needed in order to show the existence of $W^*(A)$.   Although such 
`existence proofs' are standard fare, we include
most of the details below for the readers convenience.    

Let $A$ be a dual operator algebra.  We suppose that
the cardinality of $A$ is less than or equal to
a large enough
infinite, uncountable, cardinal $I$, and define ${\mathcal F}$ to be
the set of normal completely contractive representations
$\pi : A \rightarrow B(\ell^2(J))$ where $J$ varies over 
the  cardinals corresponding to subsets of $I$.
We write $H_\pi = \ell^2(J)$.
Define $j = \oplus \{\pi : \pi \in {\mathcal F} \}$,
that is, 
$j(a) = \oplus_{\pi \in {\mathcal F}} \; \pi(a)$ for all $a \in A$.   
This is a normal completely contractive representation of $A$ on a
Hilbert space $H^w = \oplus_{\pi \in {\mathcal F}} \; H_\pi$.  
In fact $j$ is also completely
isometric, as may be seen by the standard arguments 
(Sketch: take any one normal completely isometric representation 
$\sigma$ on a Hilbert space 
$H$.  If $H$ is of dimension $\leq I$, we are done.
If not, then for each finite 
subset $F$ of $H$ set $H_F$ to be the Hilbert space 
generated by $\sigma(A) F$, and set $\pi_F$ 
to be $\sigma(\cdot)_{|_{H_F}}$.   Each 
$\pi_F$ is  unitarily equivalent to a representation on 
${\mathcal F}$.  Also, for each $x = [x_{ij}] \in M_n(A)$ we 
may clearly find such a finite set $F$ such that 
the norm of $[\pi_F(x_{ij})]$ is close to that of 
$[\sigma(x_{ij})]$.)
Thus, by the Krein-Smulian theorem,
 $j$ is a homeomorphism for the weak* topologies, with
weak* closed range.  The projection of $H^w$ onto its
`$\pi$'th coordinate will be written as $P_\pi$.
We define $W^*(A)$ to be the
von Neumann algebra inside $B(H^w)$ 
generated by $\{ j(a) : a \in A \}$.
If $\theta : A \rightarrow B(H)$ is any
normal completely contractive representation of $A$, with dimension
$H \leq I$, then there is a unitary $U$ such that
$\rho = U^* \theta(\cdot) U \in {\mathcal F}$.                              
Define $\tilde{\rho} : W^*(A) \rightarrow B(H_\rho)$ to
 be $\tilde{\rho}(T) = P_\rho T_{|_{H_\rho}}$.   Then
$\tilde{\rho}$ is a weak* continuous *-homomorphism on
$W^*(A)$, and $\tilde{\rho} \circ j = \rho$.
Then $\tilde{\theta} =
U \tilde{\rho}(\cdot) U^*$ is a weak* continuous *-homomorphism
$W^*(A) \rightarrow B(H)$, and $\tilde{\theta}  \circ j =
\theta$.  Clearly $\tilde{\theta}$ is the unique such  *-homomorphism.

Thus we have shown that $W^*(A)$ has the desired
universal
property 
at least for representations on Hilbert spaces of dimension $\leq I$.  
From this fact and a routine Zorn's lemma argument it is not hard to show 
that $W^*(A)$ has the desired universal property for 
arbitary dimensions of normal representations.  

We next prove a variant on Proposition \ref{Prop1}:

\begin{proposition}  \label{Pro1'}  Suppose that
$\rho : W^*(A) \rightarrow B(H)$ is a normal *-representation
with the property that for every normal state $\phi$ on $W^*(A)$,
there is an $x \in H$ such that $\phi = < \rho(\cdot) x , x >$
on  $W^*(A)$.
Then $\rho(A)'' \subset alg lat \rho(A)$.
\end{proposition}

\begin{proof}  Notice that
$K = [\rho(A)x]$ is a Hilbert space of cardinality
$\leq I^N = I$.
Thus $dim(K)  \leq I$, and so there is a normal *-representation
$\pi$ of $W^*(A)$ on $K$ extending $\rho(a)_{|_K}$.
The rest of the proof is the same as that of \ref{Prop1b}
combined with  \ref{Prop1}.
\end{proof}

{\bf Remark:}  If $\pi$ is any faithful normal *-representation  
of $W^*(A)$ on a Hilbert space $K$, then 
$\pi(\cdot) \otimes I_\infty$ satisfies the hypothesis 
of 
Proposition \ref{Pro1'}.

\begin{corollary} \label{Co'}  For a dual operator algebra $A$,
and for any faithful normal *-representation $\pi$
of $W^*(A)$, we have $\pi(A)'' = \pi(A)$.
\end{corollary}

\begin{proof}  This is almost identical to the proof of \ref{dct1}.
\end{proof}

One may show using standard facts (see eg. 1.3 in \cite{Ri})
that any faithful normal representation of $W^*(A)$ 
restricts to a normal $A$-universal representation of $A$.
One may also prove `normal' analogues of  
parts (1)-(3) of Theorem \ref{Thm},
for a dual operator algebra $A$ and 
normal $A$-universal representations.
Indeed  this follows exactly the  proof of Theorem \ref{Thm}.
            
\section{Complements and examples}
  
If $A$ is a $\ca$ then the $A$-universal representations
are quite well understood.  We will recap some facts almost all
of which may be found in Section 1 of \cite{Ri}, and then 
we will contrast these with 
the situation for nonselfadjoint algebras.
Recall, from Section 1 above, that a module $H \in \;
_AHMOD$ is a {\em generator} (resp.  {\em cogenerator})
for $_AHMOD$ if for every nonzero morphism $R : K \rightarrow L$ 
of $_AHMOD$, there 
exists a  morphism $T : H  \rightarrow K$
(resp. $T : L  \rightarrow H$) of $_AHMOD$ with $R T \neq 0$
(resp. $T R  \neq 0$).   We shall say that 
$H$ is a {\em semigenerator}  (resp.  {\em semicogenerator}) 
 if the condition in the last sentence is 
valid in the case that $R$ is the identity map
on a nonzero Hilbert module.   Thus,  for example, $H$ is a 
semigenerator if for every nonzero Hilbert module $K$ there 
is a nonzero bounded $A$-module map $T : H \rightarrow K$.
We shall say that $H$ is a {\em  *-generator} (resp. 
an {\em  *-semigenerator}) if the definition above for 
generator (resp. semigenerator) is modified so that 
the map $T$ considered there is required
to be adjointable (that is, $T^*$ is also an $A$-module map).
      
If $A$ is a $\ca$ then generators, cogenerators, and *-generators,
are the same thing, due to the fact that 
$T$ is a bounded $A$-module map if and only if $T^*$ is one
 also.
Similarly, for $\cas$ semigenerators, *-semigenerators, and semicogenerators,
coincide.   Indeed one has:

\begin{theorem}  \label{Ri}  (\cite{Ri} Section 1)
Let $A$ be a $C^*$-algebra, and let $\pi :
A \rightarrow B(H)$ be a nondegenerate
*-representation.
View $H$ as a Hilbert $A$-module in the usual way.
The following are equivalent:
\begin{itemize}
\item [(i)]  the canonical 
(and unique) 
weak* continuous map $\bar{\pi} :
A^{**}   \rightarrow B(H)$ extending 
$\pi$, is 1-1,
\item [(ii)]  $H$ is a semigenerator for $_AHMOD$,  
\item [(iii)]  $H$ is a generator for $_AHMOD$,   
\item [(iv)]  $H$ is $A$-universal.
\end{itemize}
\end{theorem}

The above is quite useful.  For example it follows immediately
from (i) that for a finite dimensional $\ca$,  the $A$-universal
representations are exactly the 
faithful unital *-representations.

Next, let $A$ be a nonselfadjoint operator algebra, and 
let $H$ be a Hilbert $A$-module, and let
$\pi : A \rightarrow B(H)$ be the associated representation.   
We consider the following
properties that $H$ may or may not have:
 
\begin{itemize}
\item [(DCP)]   $\pi$ has the double commutant property; that is
$\pi(A)'' = \overline{\pi(A)}^{weak^*}$.  
\item [(I)]  $\pi$ has the double commutant property and 
the canonical (and unique) weak* continuous map $\bar{\pi} :
A^{**}   \rightarrow B(H)$ extending
$\pi$, is completely isometric.   
\item [(II)]   $H$ is a semigenerator for $_AHMOD$.
\item [(II)']  $H$ is a semicogenerator for $_AHMOD$.
\item [(II)'']   $H$ is a *-semigenerator for $_AHMOD$. 
  \item [(III)]   $H$ is a generator for $_AHMOD$. 
\item [(III)']   $H$ is a cogenerator for $_AHMOD$.
\item [(III)'']  $H$ is a *-generator for $_AHMOD$.
 \item [(IV)]   $H$ is $A$-universal.
\end{itemize}

The following table summarizes several
earlier observations.  We leave omitted details to the reader.  

$$\text{(IV)} \Rightarrow \text{(III)}  \Rightarrow 
\text{(II)}$$
$$ \text{(IV)} \Rightarrow \text{(III)'}
  \Rightarrow
\text{(II)'}$$
$$ \text{(IV)} \Rightarrow \text{(I)} \Rightarrow \text{(DCP)}$$
$$ \text{(IV)} \Leftrightarrow \text{(III)''} \Leftrightarrow \text{(II)''}$$
$$ \left( \text{(III)} \; \text{ or} \; \text{
 (III)'} \right)  \Rightarrow  \text{(DCP)}$$ 

\vspace{3 mm}

{\bf Example 1.}  Let 
$A = {\mathcal T}_2$ be the algebra of $2 \times 2$ upper triangular 
matrices.  In this case it is possible to precisely characterize
the representations with the double commutant property.
First notice that the representations $\pi$ of ${\mathcal T}_2$ 
on a Hilbert space $H$ are of one of the following types:
\begin{itemize}
\item [(a)]  $H = H_1 \oplus H_2$ with 
$H_1 \neq \{ 0 \}$, $H_2 \neq \{ 0 \}$ and 
$H_1 \perp H_2$; and there 
exists a contraction $T : H_2 \rightarrow H_1$ such that 
$\pi(A)(\zeta + \eta) = a_{11} \zeta + a_{12} T(\eta) +
a_{22} \eta$, for all $\zeta \in H_1, \eta \in H_2$,
and $A = [a_{ij}] \in {\mathcal T}_2$.  
\item [(b)]  $H \neq \{ 0 \}$ and 
$\pi(A)(\zeta) = a_{11} \zeta$, for 
all $\zeta \in H$ and $A = [a_{ij}] \in {\mathcal T}_2$.
\item [(c)]  $H \neq \{ 0 \}$ and $\pi(A)(\zeta) = a_{22} \zeta$,
 for all $\zeta \in H$ and $A = [a_{ij}] \in {\mathcal T}_2$.
\item [(d)]  $H = \{ 0 \}$.
\end{itemize}   

We will not discuss the trivial case (d) below.  Clearly 
types (b) and (c) possess the double commutant property.
We write a representation $\pi$ as in (a) above
as a 3-tuple $(H_1,H_2,T)$.   

\begin{proposition} \label{uptrp}   Let 
$\pi$ be a type (a) representation of ${\mathcal T}_2$, 
with associated 3-tuple $(H_1,H_2,T)$ as above. 
\begin{itemize}
\item [(1)]  $\pi$ possesses the
double commutant property if and only if $T :
H_2 \rightarrow H_1$ is {\em not} invertible.  
\item [(2)]  $\pi$ is a semigenerator if and only if
$T(H_2)$ is not dense in $H_1$.
\item [(3)]  $\pi$ is a semicogenerator if and only if
$T$ is not 1-1.
\item [(4)]  $\pi$ is a generator if and only if
$T(H_2)$ is not dense in $H_1$ and $T \neq 0$.
\item [(5)]  $\pi$ is a cogenerator if and only if
$T$ is not 1-1 and $T \neq 0$.
\item [(6)]   $\pi$ is (completely) sub-tracing if and only if
$T(H_2)$ is not dense in $H_1$.   
\end{itemize}
\end{proposition} 

\begin{proof}  (1).
An elementary computation shows that  $\pi(A)'$ consists of
all operators $A \oplus D$, where $A \in B(H_1), D \in B(H_2)$,
such that $A T = T D$.   One observation which will be useful
later is that if $\zeta \in H_2, \eta \in H_1$ then
$A = T \zeta \otimes \eta$ and $D = \zeta \otimes T^* \eta$
satisfies $A T = T D$.  

If $T$ is invertible then the set of solutions $(A,D)$ to the 
equation $A T = T D$ is  $\{ (A,T^{-1} A T) : A \in B(H_1) \}$.
In this case the operator $z$ defined to be 
$T^{-1}$ on $H_1$ and zero on $H_2$, is easily seen to be 
in $\pi(A)''$, since $z(A \oplus T^{-1} A T)(\zeta + \eta)
= T^{-1} A \zeta = (A \oplus T^{-1} A T) z (\zeta + \eta)$
for $\zeta \in H_1, \eta \in H_2$. 
However $z$ is clearly
not in $\overline{\pi(A)}^{weak*} = \pi(A)$.

On the other hand, suppose that  $T$ is not invertible.
If $T$ is the zero operator, then any $A \in B(H_1), D \in B(H_2)$
satisfies  $A T = T D$, from which it is easily seen
that $\pi(A)'' = \pi(A)$.   Thus we may suppose that 
$T \neq 0$.    
An operator $R$ in $\pi(A)''$ may
be written as a $2 \times 2$ operator matrix with respect to
the decomposition $H_1 \oplus H_2$.  The 1-1 entry $x$  
of this matrix must therefore commute with any $A$ satisfying 
$A T = T D$ as above.   Picking 
$A = T \zeta \otimes \eta$ and $D = \zeta \otimes T^* \eta$
as above, we have that 
$x T \zeta \otimes \eta = T \zeta \otimes x^* \eta$,
so that $\eta \otimes x T \zeta  = 
x^* \eta \otimes T \zeta$.  It follows from this
that $x^* \eta$ is a scalar multiple of $\eta$ for every vector 
$\eta$, which implies that $x \in \mathbb{C} I_{H_1}$.
A similar argument shows that the 2-2 entry 
of $R$ is in $\mathbb{C} I_{H_2}$.  
A similar argument shows that the 1-2 entry
$y$ of $R$ is a scalar multiple of $T$.
To complete the proof, we need to show that the 
2-1 entry $z$ of $R$ is zero.  
The fact that $D z = z A$ yields as above that  
$\zeta \otimes z^* T^* \eta = z T \zeta \otimes \eta$
for all $\zeta \in H_2, \eta \in H_1$ as above.
It follows that either $z T = T z = 0$, or  that $T$ is
both left and right invertible.  The latter is impossible,
by hypothesis.  Thus if $T(H_2)$ is dense in $H_1$
then  $z = 0$.  On the other hand, if
$T(H_2)$ is not dense in $H_1$, then set $D = 0$ 
and let $A = \xi \otimes \sigma$, where $\sigma
\in T(H_2)^\perp$.   Clearly $A T = 0 = T D$, so that
$z A = z \xi \otimes \sigma = 0$.  Since $\xi$ is
arbitrary we must have $z = 0$.    

We leave (2)-(6) as simple but tedious exercises.  For 
example, to check (4) one assumes that $T(H_2)$ is not 
dense in $H_1$, and then one considers the 
various cases 
that can arise (corresponding to the 
types (a)-(c) of ${\mathcal T}_2$-modules) for nonzero maps $R$ between 
Hilbert ${\mathcal T}_2$-modules.  
\end{proof}

From the above it is easy to find very simple finite dimensional 
completely isometric representations of ${\mathcal T}_2$
satisfying (I), but not (II), (II)', (III), (III)', or (IV).  
Similarly (II), or even (II) together with (II)', 
does not imply (IV).   And (III), or even (III) together with (III)',
does not imply (IV). 
Indeed by the last result one can easily 
find finite dimensional completely 
isometric representations satisfying 
(I), (III), and (III)'; however
no finite dimensional representation of 
${\mathcal T}_2$ can be $A$-universal.  To see this
note that by (1) of  
\ref{uni} any $A$-universal representation on a
Hilbert space $H$ is
$C^*(A)$-universal, and therefore extends 
to a faithful representation of $C^*(A)$ on $H$.  Indeed 
the $A$-universal representations are,
by (1) of \ref{uni} and \ref{Ri}, in 1-1 correspondence 
with the normal faithful *-representations of 
$C^*(A)^{**}$.
However in \cite{BOMD} Section 2 it is shown that 
$C^*({\mathcal T}_2)$ is infinite dimensional.
 
\vspace{3 mm}
                 
{\bf Example 2.}    We consider a generalization of 
Example 1, which will show for example that 
(II) does not imply the double commutant property,
and also that even for completely isometric representations
(II) and (III)
 may differ.

Let $X$ be an operator space, and let ${\mathcal U}(X)$ 
be the canonical `upper triangular' algebra consisting 
of upper triangular $2 \times 2$ matrices with 
scalars on the diagonal and $X$ in the 1-2 corner.  
Then ${\mathcal U}(X)$ has a canonical operator space 
structure making it a unital operator algebra.  See the 
last section in \cite{BOMD} for example.   As is
spelled out there,  the nontrivial representations $\pi$ of 
${\mathcal U}(X)$ are in 1-1 correspondence with 
completely contractive maps $\alpha : X 
\rightarrow B(H_2,H_1)$.   Also, 
$\pi$ is completely isometric if and only if 
 $\alpha$ is completely isometric.  We may thus 
associate with $\pi$ the tuple $(H_1,H_2,\alpha)$.
If $X = \mathbb{C}$ then 
this is simply the $(H_1,H_2,T)$ notation we met in 
Example 1.    If $\pi$ is such a representation, with 
$H_1 \neq \{ 0 \}$ and $H_2 \neq \{ 0 \}$, then 
it is easy to compute the commutant $\pi(A)'$.  Analoguously
to Example 1, this commutant consists of the 
operators $A \oplus D$ with $A \in B(H_1), D  \in B(H_2)$,
such that $A \alpha(x) = \alpha(x) D$ for all $x \in X$.
The second commutant $\pi(A)''$ therefore is the 
set of $2 \times 2$ operator matrices 
$$\left[ \begin{array}{cl} x & y \\ z & w \end{array} \right]$$
satisfying the equations 
$Ax = xA$,  $Dw = wD, Ay = yD$ and $zA = Dz$,
 whenever $A \in B(H_1), D  \in B(H_2)$ satisfy 
$A \alpha(x) = \alpha(x) D$ for all $x \in X$.  From this
and Theorem \ref{Thm}
 (resp. \ref{Co'}) we can deduce `double
commutant theorems' for $X$.  For example it follows that:

\begin{corollary} \label{Thm3}  For any operator space 
(resp. dual operator space) $X$, there exists a completely isometric 
(resp. and weak* homeomorphic) linear $\alpha : X
\rightarrow B(H_2,H_1)$ such that the weak* closure of 
$\alpha(X)$ (resp. such that $\alpha(X)$) coincides with the set of operators 
$S \in B(H_2,H_1)$ such that 
$A S = S D$, whenever  $A \in B(H_1), D  \in B(H_2)$ satisfies
$A T = T D$ for all $T \in \alpha(X)$.    
\end{corollary}   

It is elementary to check that a 
representation $\pi$ of ${\mathcal U}(X)$, associated with 
a tuple $(H_1,H_2,\alpha)$ as above (with $H_1$ and 
$H_2$ nonzero), is a semigenerator 
if and only if the span of the ranges of the 
operators $\alpha(x)$, for all $x \in X$,
is not dense in $H_1$.    Also,
the representation is a semicogenerator
if and only if $\cap_{x \in X} \; \text{Ker} \; \alpha(x)
\; \neq \{ 0 \}$.   From this it is quite easy to find semigenerators
or semicogenerators which do not satisfy (DCP).  For 
example, choose 
$\alpha$ with the span of the ranges of the
operators $\alpha(x)$ not dense in
$H_1$, but for which there exist operators $w \in B(H_2)$ which are 
not scalar multiples of the identity such that 
$Dw = wD$ for all $A \in B(H_1), D  \in B(H_2)$ satisfying
$A \alpha(x) = \alpha(x) D$ for all $x \in X$.  (For a
concrete such example let $X = \ell^\infty$ and $\alpha(x) = 
S \; diag \{ x \}$, where $S$ is the forward shift.)   Then the 
$2 \times 2$ operator matrix  with $w$ in the 
2-2 entry and other entries zero, is in $\pi(A)''$ but not
in the weak* closure of  $\pi(A)$.    Thus 
(II) does not imply the (DCP).    Similar considerations
show that (II)' does not imply the (DCP).    
 
We now exhibit an example of a
completely isometric representation $\pi$ of the  
type considered in Example 2, which satisfies (II) and (II)',
but which is not a generator (i.e. does not have (III)).   
One such is given by 
$\alpha : \ell^\infty \rightarrow B(\ell_2)$ 
of the form $\alpha(x) = S \; diag \{ 0, x \}$, where 
$S$ is the forward shift again.   It is 
clear that this has property (II) and (II)'.  To see 
that (III) fails we appeal to the following Claim:
Let $X$ be a non-reflexive dual operator space.
Then any representation $\pi :
{\mathcal U}(X) \rightarrow B(H)$
associated with a 3-tuple $(H_1,H_2,\alpha)$ 
with $\alpha$  weak* continuous, is a generator.
To prove this Claim, 
consider a fixed non weak* continuous contractive linear 
functional $\beta$ on $X$,  and consider the representation 
of ${\mathcal U}(X)$ on $K = \mathbb{C}^2$ associated with the 
tuple $(\mathbb{C},\mathbb{C},\beta)$.   Let 
$L = \mathbb{C}$, with the `type (c) action' of 
 ${\mathcal U}(X)$ described in
example 1, and let $R \in B(K,L)$ 
be the projection onto the second coordinate.  This is clearly 
a nonzero $A$-module map on $K$.    A nonzero 
$A$-module map $T : H \rightarrow K$  is easily seen 
to be necessarily of the 
form $T_1 \oplus T_2$, 
where $T_1 \in B(H_1,K_1), T_2 \in B(H_2,K_2)$  may be any pair 
satisfying $T_1 \alpha(x) = 
\beta(x) T_2$ for all $x \in X$.
In fact this is true for 
any ${\mathcal U}(X)$-modules $H, K$.    If $T_2 \neq 0$ this
implies that $\beta(x)$ is a constant multiplied by 
$\langle \alpha(x) \xi , \sigma \rangle$ for some 
vectors $\xi , \sigma$.   This implies the contradiction
that $\beta$ is weak* continuous.  Thus $T_2 = 0$, so that 
$R T = 0$.

\vspace{3 mm}

The examples above rule out most of the 
variants for nonselfadjoint algebras of 
the remaining implications of \ref{Ri}.  Some
 questions 
which we have not taken the time to settle, are the following:
Does a 
completely isometric representation satisfying ((II) and (II)') 
automatically possess the double commutant property?  
Also, for faithful representations of 
unital operator algebras $A$, how close is the condition 
(DCP) to the condition ((III) or (III)')?   To 
the condition ((II)  or (II)')?     

Finally, we remark that there are other 
variants on the definition of `generator', which are
situated between (III) and (IV).  In particular the class of  
Hilbert $A$-modules $H$ with the following 
property: For any other 
Hilbert $A$-module $K$ there is a cardinal $\alpha$ and a
bounded module map $T : H^{(\alpha)} \rightarrow K$ which is
surjective (resp. has dense range, is a 1-quotient map). 
We will not say anything further about these three classes 
except that they contain (but are not equal to) the class (IV), 
and are contained in class (III), and hence are faithful 
and satisfy the double commutant property.

\section{Nonunital operator algebras}
 
In this section we verify that all of 
Theorem  \ref{Thm}, with the exception of
part (4), is valid more generally for
operator algebras $A$ which do not have a c.a.i..
We shall see that if part (4) was valid too then
$A$ must have a c.a.i..
 
We say that a
homomorphism $\pi : A \rightarrow B(H)$ is
{\em nondegenerate} if the span of terms of the form
$c_1 c_2 \cdots c_n \zeta$, for $\zeta \in H$ and
$c_i \in \pi(A) \cup \pi(A)^*$, is dense in $H$.
  Perhaps a better name for this is {\em *-nondegeneracy}, but for
simplicity we will use the other name here. 
We will not use this
fact, but
any contractive homomorphism $\pi : A \rightarrow B(K)$ can
be replaced by a nondegenerate one, by
restricting each $\pi(a)$ to the closed subspace of
$K$ densely spanned by the products $c_1 c_2 \cdots c_n \zeta$
mentioned above.
We remark that if $A$ has a c.a.i. then this new definition of
nondegeneracy of representations coincides with the old.  To see this,
suppose  that  $\{ e_\alpha \}$ is a c.a.i. for $A$, and that
$\pi$ is a contractive homomorphism which is
nondegenerate in the new sense above.  By 
fact (\ref{anc}) from Section 1 we know
that  
$\{ \pi(e_\alpha) \}$ is a c.a.i. for the $C^*$-subalgebra
of $B(H)$ generated by $\pi(A)$.
Thus $\pi(e_\alpha) \rightarrow Id$ strongly on $H$.  The
converse is easier.
 
If $A$ is any operator algebra then a recent paper
\cite{ram} proves the remarkable results that a) there
are unique matrix norms on $A^+ = A \oplus \mathbb{C}$ such that
$A^+$ is a unital abstract operator algebra (with identity
$1_+ = (0,1)$) containing $A$ completely
isometrically, and b)  given a contractive (resp. completely
contractive, isometric, completely isometric)
homomorphism $\varphi : A \rightarrow B$ between operator algebras,
the extension $\varphi^+ : A^+ \rightarrow B^+$ given by
$\varphi^+(a + \lambda 1_+)  = \varphi(a)  + \lambda 1_+$,
for $a \in A, \lambda \in \mathbb{C}$,
is also a contractive (resp.  completely
contractive, isometric, completely isometric) homomorphism.
From this it is easy to define a $C^*$-envelope (in the spirit 
of Arveson and Hamana \cite{Ham,Arv1}) and a 
maximal universal $C^*$-algebra of operator algebras without
a c.a.i..   Again for specificity
will do this in the operator space framework, as opposed to
the Banach algebra version.
 
If  $A$ is any operator algebra then we define $C^*_e(A)$ (resp. 
$C^*(A)$) to be
the $C^*$-subalgebra of $C^*_e(A^+)$ (resp. $C^*(A^+)$)
generated by the copy of $A$.   See \cite{Ham,Arv1} 
for the basic properties of the $C^*$-envelope 
$C^*_e(A^+)$.   We claim that
$C^*_e(A)$ (resp. $C^*(A)$) has the appropriate
universal properties, analogous to the
well known properties they have in the case that $A$ has a
c.a.i..  We first treat $C^*_e(A)$.   If $\pi : A \rightarrow B$
is a completely isometric homomorphism into a $C^*$-algebra $B$
such that $\pi(A)$ generates $B$ as a $C^*$-algebra,
then $\pi^+ : A^+ \rightarrow B^+$
is a completely isometric homomorphism into a $C^*$-algebra,
whose range generates $B^+$ as a $C^*$-algebra.  Thus by
the universal property of $C^*_e(A^+)$,
there is a
 surjective *-homomorphism $ \rho : B^+ \rightarrow C^*_e(A^+)$
such that $\rho \circ \pi^+ = j$, where
$j : A^+ \rightarrow C^*_e(A^+)$ is the canonical
embedding.  Let $\theta$ be $ \rho$ restricted to
$B$, then $\theta$ is a  *-homomorphism with
$$\theta(\pi(a)) = \rho(\pi^+(a)) = j(a) \in C^*_e(A)$$
for all $a \in A$.  Thus
$\theta$ maps $B$ into $C^*_e(A)$, and the above shows that
$(C^*_e(A),j)$ has the universal property which one would desire
for a `$C^*$-envelope of $A$'.    
 
We now check  that
$C^*(A)$ has the universal property which one would desire.
Suppose that 
$\pi$  is a completely contractive
homomorphism from $A$ into a $C^*$-algebra $B$.  By
the universal property of $C^*(A^+)$, there is a
*-homomorphism $\rho : C^*(A^+) \rightarrow B^+$
such that $\rho \circ \kappa = \pi^+$, where
$\kappa : A^+ \rightarrow C^*(A^+)$ is the canonical
embedding.   
Let $\theta$ be $\rho$ restricted to
$C^*(A)$; then $\theta$ is a  *-homomorphism 
with $$\theta(\kappa(a)) = \rho(\kappa(a)) = \pi^+(a) = \pi(a) \in B,$$ 
for all $a \in A$.   Thus
$\theta$ maps $C^*(A)$ into $B$.  
 
If $A$ is a nonunital operator algebra then we let $_AHMOD$ be the
category of nondegenerate Hilbert $A$-modules, using the
definition of `nondegenerate' given at the beginning of this section.
By the universal property of $C^*(A)$, the objects in $_AHMOD$ 
are `the same as' the objects in $_{C^*(A)}HMOD$.   In particular, a
completely contractive 
representation of $A$ is nondegenerate in the new 
sense if and only if the associated representation of  $C^*(A)$
is nondegenerate in the usual sense.  
We may define
direct sums in $_AHMOD$ by associating them with the
corresponding direct sums in $_{C^*(A)}HMOD$.   
Thus, a direct sum
of Hilbert $A$-modules is nondegenerate if and only if
every one of the individual summand Hilbert $A$-modules is nondegenerate.
The fact from Section 3 labelled 
(\ref{adj}) is still valid with the same proof, and so we may treat
$A$-complemented submodules and direct summands in $_AHMOD$
just as we did before.  Indeed Corollary \ref{uni} also
carries over verbatim, as does (1)-(3) of Theorem \ref{Thm}.
We define the {\em universal representation} $\pi_u$ of
$A$ to be the restriction to $A$ of the universal representation
$\pi_u$ of $C^*(A)$.
The facts in the second paragraph of Section 4 also transfer
immediately, the only difference being that $A^{**}$ and
$\theta$ there need not be unital.  Now we see that (5) of
Theorem \ref{Thm} carries over verbatim too.   Thus all
of Theorem \ref{Thm}, with the exception of
part (4), is valid when $A$ is an  operator algebra with no c.a.i..
 
Indeed it is clear that if (4) and (5)
of  \ref{Thm} both hold, then $A^{**}$ is unital, and from the
theory of Banach algebras 
it follows that $A$ has a c.a.i..   
 
\begin{corollary}
\label{chnu}  An operator algebra $A$ possesses a c.a.i. if
and only if  
for every nondegenerate contractive homomorphism
 $\pi : A \rightarrow B(H)$,
we have $x \in [\pi(A) x]$ whenever $x \in H$.
\end{corollary}

\begin{proof}     The one direction is easy, using 
the facts noted at the beginning of this section.
The other direction may be proved by noting that if the
hypothesis holds then the rest of the proof of Theorem \ref{Thm}
(4)
is easily amended to yield the nonunital case.  Hence $A^{**}
\cong A''$ is unital, and so 
$A$ has a c.a.i. 
as mentioned above \ref{chnu}.  \end{proof}

{\bf Remarks:}  A direct `reflexivity' proof  of \ref{chnu} may also 
be given.    Also, we point out that
the qualification `for all' in \ref{chnu} may not be replaced
by `for some completely isometric nondegenerate
homomorphism $\pi : A \rightarrow B(H)$'.  To see this consider
the unitary operator $U$ on $H = L^2[0,2\pi]$ given
by an irrational rotation.  If $A \subset B(H)$ is the
uniform closure of the span of $U, U^2, U^3, \cdots$, then it
is fairly clear that $f \in [A f]$ for all $f \in L^2[0,2\pi]$.
In particular $A$ acts nondegenerately on $H$.
If $A$ contained a c.a.i. $\{ E_\alpha \}$ then
$E_\alpha U \rightarrow U$, and so $
E_\alpha \rightarrow Id_H$.  Hence $Id_H \in A$.
On the other hand,  $A$ is contained in the closure of the $C^*$-algebra
generated by $I$ and $U$, and this latter $C^*$-algebra
is isomorphic to the set $C(\mathbb{T})$
of continuous functions on the
circle, by basic spectral theory.   Under this
isomorphism $A$ corresponds to the {\em nonunital} ideal $z A(\mathbb{D})$ in
the disk algebra $A(\mathbb{D})
\subset C(\mathbb{T})$.   Here $z$ represents the function
$e^{i \theta } \mapsto e^{i \theta }$ on $\mathbb{T}$.   This
contradiction shows that $A$ does not have a c.a.i..

\vspace{3 mm}

Acknowledgments:   This work was begun during the 2000 summer 
workshop at Texas A \& M University, and
we thank the organizers of that program for their support
and venue.   We acknowledge the influence of 
\cite{Ri}.   We thank C. Le Merdy and
P. S. Muhly 
for several helpful inputs. 
Indeed the possibility of results
such as these was discussed with Le Merdy in 1999.  Also, he has pointed 
out to us the work \cite{Gi} of J. A. Gifford in which a double commutant 
theorem is obtained for representations satisfying a certain 
condition.  This condition has no real overlap with ours, and it 
seems likely that an operator algebra would need to be very 
close to being a $C^*$-algebra (eg. similar to
a $C^*$-algebra) in order to satisfy his condition.   Indeed 
Gifford's intention was to prove this latter assertion, which 
he does in certain cases.   
Our theorem
\ref{pThm} does imply Gifford's double commutant result.    
It also implies that certain induced representations
of tensor algebras studied by the second author with P. S. 
Muhly, have the double commutant property.
We have not tried to connect our ideas
with the asymptotic double commutant theorems for
$C^*$-algebras of
Hadwin (see eg. \cite{Had}).
There is possibly also a connection to
Hadwin's operator
valued spectrum.

\end{document}